\documentstyle[11pt,leqno]{article}

\newtheorem{theorem}{{\sc Theorem}}

\newcommand{\bt}{\begin{theorem}}

\newcommand{\et}{\end{theorem}}

\newcommand{\newsection}[1]{\setcounter{equation}{0} \setcounter{theorem}{0}

\section{#1}}

\newcommand{\NI}{\noindent}

\newcommand{\bea}{\begin{eqnarray}}

\newcommand{\eea}{\end{eqnarray}}

\def \spec#1 {\mathop{#1}}

\def \b #1 {\bf #1}

\newcommand {\nnb }{\nonumber}

\newcommand {\CC}{\centerline}

\newcommand{\cla}{{\cal A}}

\newcommand{\clf}{{\cal F}}

\newcommand{\cll}{{\cal L}}

\newcommand{\clu}{{\cal U}}

\newcommand{\clw}{{\cal W}}

\newcommand{\cly}{{\cal Y}}

\newcommand{\clz}{{\cal Z}}
\newcommand{\ve} {\varepsilon}

\newcommand{\bta}{\beta}

\newcommand{\pal}{\partial}

\newcommand{\Omg}{\Omega}

\newcommand{\lmd}{\lambda}

\newcommand{\tha}{\theta}

\newcommand{\Tha}{\Theta}

\newcommand{\al}{\alpha}

\newcommand{\ity}{\infty}

\newcommand{\raro}{\rightarrow}

\newcommand{\vsp}{\vskip 1em}

\newcommand{\be}{\begin{equation}}

\newcommand{\ee}{\end{equation}}

\newcommand{\ben}{\begin{eqnarray*}}

\newcommand{\een}{\end{eqnarray*}}

\begin{document}
\CC {\bf Parametric Estimation for Processes Driven by }
\CC{\bf Infinite Dimensional Mixed Fractional Brownian Motion}
\vsp
\CC{\bf B.L.S. Prakasa Rao}
\vsp
\CC{\bf CR Rao Advanced Institute of Mathematics, Statistics}
\CC{\bf and Computer Science, Hyderabad 500046, India}
\vsp
\NI{\bf Abstract:} Parametric and nonparametric inference for stochastic processes driven by a fractional Brownian motion were investigated in Mishura (2008) and Prakasa Rao (2010) among others. Similar problems for processes driven by an infinite dimensional fractional Brownian motion were studied in Prakasa Rao (2004, 2013), Cialenco et al. (2009) and others. Parametric estimation for processes driven by an infinite dimensional mixed fractional Brownian motion is discussed in this article.  
\vsp
\newsection{Introduction}
Statistical inference  for diffusion type processes satisfying
stochastic differential equations driven by Wiener processes has
been studied earlier and a comprehensive survey of various methods
is given in Prakasa Rao (1999a) . There has been a recent interest
to study similar problems for stochastic processes driven by a
fractional Brownian motion to model processes having long range dependence. Le Breton (1998) studied parameter
estimation and filtering in a simple linear model driven by a
fractional Brownian motion. Kleptsyna and Le
Breton (2002) studied parameter estimation problems for fractional
Ornstein-Uhlenbeck type process driven by a fractional Brownian motion. This is a fractional analogue of
the Ornstein-Uhlenbeck process driven by a standard Wiener process. It is a continuous time first
order auto-regressive process $X=\{X_t, t \geq 0\}$ which is the
solution of a one-dimensional homogeneous linear stochastic
differential equation driven by a fractional Brownian motion (fBm)
$W^H= \{W_t^H, t \geq 0\}$ with Hurst parameter $H \in [1/2, 1).$
Such a process is the unique Gaussian process satisfying the
linear integral equation \be X_t= \theta \int_0^tX_s ds + \sigma
W_t^H, t \geq 0. \ee They investigated the problem of estimation of
the parameters $\theta$ and $\sigma^2$ based on the observation
$\{X_s, 0 \leq s \leq T\}$ and proved that the maximum likelihood
estimator $\hat \theta_T$ is strongly consistent as $T \raro
\ity.$ More general classes of stochastic processes
satisfying linear stochastic differential equations driven by a
fractional Brownian motion were studied and the asymptotic properties of
the maximum likelihood and the Bayes estimators for parameters
involved in such processes is investigated in Prakasa Rao (2003). Prakasa Rao (2010) gives a comprehensive
discussion on problems of estimation for processes driven by a fractional Brownian motion.

Geometric Brownian motion driven by a standard Brownian motion has been widely used for modeling fluctuations of share prices in a stock market using Black-Scholes model. However efforts to model fluctuations in financial  markets with long range dependence through processes driven by a fractional Brownian motion were not successful as it was noted that such a modeling creates arbitrage opportunities contrary to the fundamental assumption of no arbitrage opportunity for modeling rational market behaviour. Cheridito (2001) proposed modeling through processes driven by a mixed fractional Brownian motion. It was shown by Cheridito (2001) that a mixed fractional Brownian motion is a semimartingale if and only if the Hurst index $H$ is either equal to $\frac{1}{2}$ reducing the process to a Wiener process or $H\in (3/4,1)$. Furthermore the probability measure generated by such a process is absolutely continuous with respect to the probability measure generated by a Wiener process if $H=1/2$ or $H\in (3/4,1).$ This in turn will lead to no arbitrage opportunities for modeling financial market behaviour through processes driven by a mixed fractional Brownian motion. This discussion is to motivate the study of processes driven by a mixed fractional Brownian motion.  

The problem of estimation of parameters for processes driven by processes which are mixtures of  independent Brownian and fractional Brownian motions started from the works of Cheridito (2001), Rudomino-Dusyatska (2003) and more recently in Prakasa Rao (2015a,b;2017a,b; 2018a,b; 2019, 2020, 2021a,b) among others. Mixed fractional Brownian models were  studied in Mishura (2008) and Prakasa Rao (2010) . Cai et al. (2016) present a new approach via filtering for analysis of mixed processes of type $\{X_t=B_t+G_t, 0\leq t \leq T\}$ where $\{B_t, 0\leq t \leq T\}$ is a Brownian motion and $\{G_t, 0\leq t \leq T\}$ is an independent Gaussian process. Statistical analysis of mixed fractional Ornstein-Uhlenbeck process was investigated in Chigansky and Kleptsyna (2019). Fractional Ornstein-Uhlenbeck type process driven a mixed fractional Brownian motion  has also been termed as ``mixed fractional Ornstein-Uhlenbeck process"' in Marushkevych (2016). Large deviations for drift parameter estimator of a mixed fractional Ornstein-Uhlenbeck process were studied by Marushkevych (2016).

Huebner et al. (1995) initiated the study of parametric estimation for a class of stochastic partial differential equations (SPDE) in the presence of white noise or the driving force is an infinite dimensional Wiener process. These results were extended to parabolic stoochastic partial differential equations in Huebner and Rozovskii (1995). Prakasa Rao (2000) studied Bayes estimation for stochastic partial differential equations in the white noise case. For other results on parametric inference for SPDEs, see Prakasa Rao (2000, 2001, 2002a,b,2004, 2013). A comprehensive survey of results is given in Prakasa Rao (2001,2002). Parameter estimation for a two-dimensional  stochastic Navier-Stokes equation driven by infinite dimensional fractional Brownian motion was studied in Prakasa Rao (2013). Lototsky and Rozovsky (2017) give an extensive survey of theory of SPDE and a discussion on parametric inference for such processes. Cialenco and his coworkers obtained several results dealing with parametric inference for SPDE based on continuous observation or discrete sampling of the processes. Cialenco (2019) gives a survey of their results. 

Our aim in this paper is to study parametric inference for processes driven by infinite dimensional mixed fractional Brownian motion. As far as we are aware, this problem has not been investigated earlier. 

\vsp
\newsection{Properties of processes driven by a mfBm}
Let $(\Omega, \clf, (\clf_t), P) $ be a stochastic basis
satisfying the usual conditions. The natural filtration of a
stochastic process is understood as the $P$-completion of the
filtration generated by this process. Let $\{W_t, t \geq 0\}$ be a standard Wiener process and $W^H= \{W_t^H, t \geq 0 \}$ be an independent  normalized  fractional Brownian motion with Hurst parameter $H \in (0,1)$, that is, a Gaussian process with continuous sample paths such that $W_0^H=0, E(W_t^H)=0$ and
\be
E(W_s^H W_t^H)= \frac{1}{2}[s^{2H}+t^{2H}-|s-t|^{2H}], t \geq 0, s \geq 0.
\ee
Let
$$\tilde W_t^H= W_t+ W_t^H, t \geq 0.$$
The process $\{\tilde W_t^H, t \geq 0\}$ is called the {\it mixed fractional Brownian motion} with Hurst index $H.$ We assume here after that Hurst index $H$ is {\it known} and that $H \in (\frac{3}{4},1).$ 
\vsp
Let us consider a stochastic process $Y=\{Y_t, t \geq 0\}$ defined
by the stochastic integral equation \be Y_t= \int_0^t C(s) ds  +\tilde W_t^H, t \geq 0 \ee where the process $C=\{C(t), t
\geq 0\}$ is an $(\clf_t)$-adapted process. For convenience, we write
the above integral equation in the form of a stochastic
differential equation
\be
dY_t= C(t) dt + d\tilde W_t^H, t \geq 0
\ee
driven by the mixed fractional Brownian motion $\tilde W^H.$ Following the recent works by Cai et al. (2016) and Chigansky and Kleptsyna (2019), one can construct an integral transformation that transforms the mixed fractional Brownian motion $\tilde W^H$ into a martingale $M^H.$ Let $g_H(s,t)$ be the solution of the integro-differential equation
\be
g_H(s,t)+H \frac{d}{ds}\int_0^t g_H(r,t)|s-r|^{2H-1} sign(s-r)dr=1, 0<s<t.
\ee
Cai et al. (2016) proved that the process
\be
M_t^H= \int_0^tg_H(s,t)d\tilde W_s^H, t \geq 0
\ee
is a Gaussian martingale with quadratic variation
\be
<M^H>_t= \int_0^tg_H(s,t)ds, t \geq 0
\ee
Let $w^H_t$ denote the quadratic variation $<M^H>_t$ over the interval $[0,t].$ It is known that the natural filtration of the martingale $M^H$ coincides with that of the mixed fractional Brownian motion $\tilde W^H.$ Suppose that, for the martingale $M^H$ defined by the equation (2.5), the sample paths of the process $\{C(t), t \geq 0\}$ are smooth enough in the sense that the process
\be
Q_H(t)= \frac{d}{d<M^H>_t}\int_0^tg_H(s,t)C(s)ds, t \geq 0
\ee
is well defined. Define the process
\be
Z_t= \int_0^tg_H(s,t) dY_s, t \geq 0.
\ee
As a consequence of the results in Cai et al. (2016), it follows that the process $Z$ is a fundamental semimartingale associated with the process $Y$ in the following sense.
\vsp
\NI{\bf Theorem 2.1:} {\it Let $g_H(s,t)$ be the solution of the equation (2.4). Define the process $Z$ as given in the equation (2.8). Then the following relations hold.

\NI(i) The process $Z$ is a semimartingale with the decomposition
\be
Z_t= \int_0^tQ_H(t)d<M^H>_s + M^H_t, t \geq 0
\ee
where $M^H$ is the martingale defined by the equation (2.5).

\NI(ii) The process $Y$ admits the representation
\be
Y_t=\int_0^t\hat g_H(s,t)dZ_s, t \geq 0
\ee
where
\be
\hat g_H(s,t)= 1-\frac{d}{d<M^H>_s}\int_0^tg_H(r,s)dr.
\ee
\NI (iii) The natural filtrations $(\cly_t)$ and $(\clz_t)$ of the processes $Y$ and $Z$ respectively coincide.}
\vsp
Applying Corollary 2.9 in Cai et al. (2016), it follows that the probability measures $\mu_Y$ and $\mu_{\tilde W^H}$ generated by the processes $Y$ and $\tilde W^H$ on an interval $[0,T]$ are absolutely continuous   with respect to each other and the Radon-Nikodym derivative is given by
\be
\frac{d\mu_Y}{d\mu_{\tilde W^H}}(Y)= \exp[\int_0^TQ_H(s) dZ_s-\frac{1}{2}\int_0^T[Q_H(s)]^2 d<M^H>_s]
\ee
which is also the likelihood function based on the observation $\{Y_s,0\leq s \leq T.\}$ Since the filtrations generated by the processes $Y$ and $Z$ are the same, the information contained in the families of $\sigma$-algebras  $(\cly_t)$ and $(\clz_t)$  is the same and hence the problem of the estimation of the parameters involved  based on the observation $\{Y_s, 0\leq s \leq T\}$ and $\{Z_s, 0\leq s \leq T\}$ are equivalent. 
Since the process $\{Z_s, 0\leq s \leq T\}$ is driven by a martingale, it is convenient to discuss asymptotic behaviour of the estimators through limit theorems available for martingales. This explanation motivates the study of problem of estimation through the process $Z$ instead of the original process $Y.$
\vsp
\newsection{Parametric estimation for SPDE driven by infinite dimensional mfBm}
\vsp Kallianpur and Xiong (1995) discussed the properties of
solutions of stochastic partial differential equations (SPDE)
driven by infinite dimensional fractional Brownian motion. They
indicate that SPDE's are being used for stochastic modelling, for
instance, for  the study of neuronal behaviour in neurophysiology
and in building stochastic models of turbulence. The theory of
SPDE's is investigated in Ito (1984), Rozovskii (1990) and Da
Prato and Zabczyk (1992). Huebner et al. (1993) started the
investigation of maximum likelihood estimation of parameters of
two types of SPDE's and extended their results for a class of
parabolic SPDE's in Huebner and Rozovskii (1995). Asymptotic
properties of Bayes estimators for such problems were discussed in
Prakasa Rao (2000). A short review and a comprehensive survey of
these results are given in Prakasa Rao (2001,2002). Our aim in
this section is to study the problems of parameter estimation for
some SPDE driven by an infinite dimensional mixed fractional Brownian
motion. 
\vsp 
\NI{\bf Stochastic PDE with linear drift (absolutely continuous case)} 
\vsp Let $U$ be a real separable Hilbert space
and $Q$ be a self-adjoint positive operator. Further suppose that
the operator $Q$ is nuclear. Then $Q$ admits a sequence of
eigenvalues $\{q_n, n \geq 1\}$ with $0 < q_n $ decreasing to zero
as $n \raro \ity$ and $\sum_{n=0}^{\ity}q_n < \ity.$ In addition
the corresponding eigen vectors $\{e_n, n \geq 1\}$ form an
orthonormal basis in $U.$ We define the {\it infinite dimensional
mixed fractional Brownian motion} on $U$ with covariance $Q$ as \be
\tilde \clw_Q^H(t)=\sum_{n=0}^{\ity}\sqrt {q_n} e_n \tilde W_n^H(t) \ee where
$\tilde W_n^H, n \geq 1$ are real independent mfBm's with Hurst index $H$
. Formal definition is given in the next
section. \vsp Let $U=L_2[0,1]$ and $\clw_Q^H$ be the infinite
dimensional mfBm on $U$ with the Hurst index $H$ and with the
nuclear covariance operator $Q.$
\vsp
Consider the process
$u_{\ve} (t,x), 0 \leq x \leq 1, 0 \leq t \leq T$ governed by the stochastic
partial differential equation
\be
du_{\ve} (t,x) = (\triangle u_{\ve} (t, x) + \tha u_{\ve} (t,x)) dt + \ve \;\;
d \tilde \clw_Q^H (t,x)
\ee
where $\triangle = \frac{\pal^2}{\pal x^2}.$ Suppose that $ \ve \raro 0$
and $\tha \in \Tha \subset R.$ Suppose the initial and the boundary
conditions are given by
\bea
u_{\ve} (0,x) & = & f(x), f \in L_2 [0,1] \\
u_{\ve} (t,0) &= & u_{\ve} (t,1) = 0, 0 \leq t \leq T .
\eea
Let us consider a special covariance operator $Q$ with $e_k = \sin k \pi
x, k \geq 1$ and $\lmd_k = (\pi k)^2, k \geq 1.$ Then $\{e_k\}$ is a
complete orthonormal system with the eigenvalues $q_i = (1+ \lmd_i)^{-1},i
\geq 1$ for the operator $Q$ and $Q= (I- \triangle)^{-1}$.
\vsp
Guerra and Nualart (2008) proved an existence and uniqueness theorem for solutions of multidimensional time dependent stochastic differential equations driven by a multidimensional fractional Brownian motion with Hurst index $H>\frac{1}{2}$ and a multidimensional standard Brownian motion. Similar results were obtained by Mishura and Shevchenko (2011) and da Silva and Erraoui (2018) under weaker conditions. Mishura et al. (2019) has  given sufficient conditions for the existence and uniqueness of a  mild solution $u_{\ve}(t,x)$ for stochastic differential equation driven by an infinite dimensional mfBm.  
\vsp
We assume that sufficient conditions hold so that there exists a unique
square integrable  solution $u_{\ve}(t,x)$ of (3.2) under the
conditions (3.3)-(3.4) and consider it as a formal sum \be u_{\ve} (t,x)
= \sum_{i=1}^{\ity} u_{i \ve}(t) e_i (x). \ee It can be checked
that the Fourier coefficient $u_{i \ve} (t)$ satisfies the
stochastic differential equation \be d u_{i \ve}(t) = (\tha -
\lmd_i) u_{i \ve}(t) dt + \frac{\ve}{\sqrt{\lmd_i + 1}} d\tilde W_i^H(t),
\; 0 \leq t \leq T \ee with the initial condition \be u_{i \ve}
(0) = v_i, \; v_i = \int_0^1 f(x) e_i (x) dx . \ee Let
$P_{\tha}^{(\ve)}$ be the probability measure generated by
$u_{\ve}$ when $\tha$ is the true parameter. Suppose $\tha_0$ is
the true parameter. Observe that the process $\{u_{i\ve}(t), 0
\leq t \leq T\}$ is a mixed fractional Ornstein-Uhlenbeck type process
(cf. Marushkevych (2016), Chigansky and Kleptsyna (2019), Cai et al. (2016)). \vsp Following the notation given in the previous section, define
\be
M_i^H(t) = \int_0^tg_H(s,t) d\tilde W_i^H(s), 0 \leq t \leq T,
\ee
\be
Q_{i\ve}(t)
=\frac{\sqrt{\lambda_i+1}}{\ve}\frac{d}{dw_t^H}\int_0^t g_H(s,t)u_{i\ve}(s)ds, t \in [0,T],
\ee
\be
Z_{i\ve}(t)= (\theta -\lambda_i)\int_0^t Q_{i\ve}(s)dw_s^H + M_i^H(t), 0\leq t \leq T.
\ee
Observe that $M_i^H$ is a zero mean Gaussian martingale. Furthermore, it follow that the process $\{Z_{i\ve}(t)\}$ is a semimartingale and
the natural filtrations $(\clz_{{i\ve}_t})$ and $(\clu_{{i\ve}_t})$ of the
processes $Z_{i\ve}$ and $u_{i\ve}$ respectively coincide.
Let $P_{i\tha}^{T,\ve}$ be the probability measure generated by the
process $\{u_{i \ve}(t), 0 \leq t \leq T\}$ when $\theta $ is the true
parameter. Let $\theta_0$ be the true parameter. It follows, by the
Girsanov type theorem, that
\bea
\log \frac{dP_{i\tha}^{T,\ve}}{dP_{i\tha_0}^{T,\ve}}& = &
 \frac{\lmd_i + 1}{\ve^2} [(\tha - \tha_0) \int_0^T Q_{i\ve}(t) dZ_{i\ve}
(t) \\ \nnb
& & \;\;\;\; - \frac{1}{2} \{ (\tha - \lmd_i)^2 - (\tha_0 - \lmd_i)^2\}
\int_0^T Q_{i\ve}^2(t) dw_t^H].
\eea
Let $u_{\ve}^N(t,x)$ be the projection of the solution $u_{\ve}(t,x)$ onto
the subspace spanned by the eigen vectors $\{e_i, 1 \leq i \leq N\}.$ Then
\be
u_{\ve}^N(t,x)= \sum_{i=1}^N u_{i \ve}(t) e_i(x)
\ee
From the independence of the processes $\tilde W_i^H, 1 \leq i \leq N$ and hence
of the processes $u_{i \ve}, 1 \leq i \leq N,$ it follows that the
Radon-Nikodym derivative, of the probability measure $P_{\theta}^{N,T,\ve}$
generated by the process $u_{\ve}^N, 0 \leq t \leq T$ when $\theta$ is the
true parameter with respect to the probability measure $P_{\theta_0}^{N,T,\ve}$
generated by the process $u_{\ve}^n, 0 \leq t \leq T$ when $\theta_0$ is the
true parameter, is given by

\bea
\log \frac{dP_{\tha}^{N,T,\ve}}{dP_{\tha_0}^{N,T,\ve}}(u_\ve^N) &=&
 \sum_{i=1}^N \frac{\lmd_i + 1}{\ve^2} [(\tha - \tha_0) \int_0^T Q_{i
\ve}(t) dZ_{i \ve}(t) \\ \nnb
& & \;\;\; - \frac{1}{2} \{(\tha - \lmd_i)^2 - (\tha_0 -
\lmd_i)^2\} \int_0^T Q_{i\ve}^2(t) dw_t^H].
\eea

Furthermore the Fisher information is given by
\bea
I_{N\ve}(\theta) & = & E_\theta [\frac{\partial \log
\frac{dP_{\theta}^{N,T \ve}}{dP_{\theta_0}^{N,T,\ve}}}{\partial \theta}]^2
\\ \nnb
& = & \sum_{i=1}^N \frac {\lambda_i+1}{\ve^2}E_\theta \{\int_0^T Q_{i\ve}^2(t)dw_t^H\}.
\eea

It is easy to check that the maximum likelihood estimator $\hat \theta_{N,
\ve}$ of the parameter $\theta$ based on the projection $u_\ve^N$ of
$u_\ve$ is given by
\be
\hat \theta_{N,\ve} = \frac {\sum_{i=1}^N(\lmd_i + 1)
 \int_0^T Q_{i\ve}(t) dZ_{i \ve}(t)}{\sum_{i=1}^N (\lmd_i+1)
\int_0^TQ_{i\ve}^2(t) dw_t^H}.
\ee
Suppose $\theta_0$ is the true parameter. It is easy to see that
\be
\ve^{-1}(\hat \theta_{N,\ve}-\theta_0)= \frac {\sum_{i=1}^N \sqrt{\lmd_i +
1} \int_0^T Q_{i\ve}(t) dM_{i}^H(t)}{\sum_{i=1}^N (\lmd_i+1)
\int_0^TQ_{i\ve}^2(t) dw_t^H}.
\ee

Observe that $M_i, 1 \leq i \leq N$ are independent zero mean Gaussian
martingales with $<M_i>=w^H, 1 \leq i \leq N.$
\vsp
\NI {\bf Theorem 3.1 :} The maximum likelihood estimator $\hat \theta_{N,\ve}$
is strongly consistent, that is,
\be
\hat \theta_{N,\ve} \raro \theta_0 \;\;\mbox{a.s}\;\; [P_{\theta_0}] \;\mbox
{as}\;\; \ve \raro 0
\ee
provided
\be
\sum_{i=1}^N \int_0^T (\lambda_i+1)Q_{i \ve}^2(t) dw_t^H \raro \ity\;\; \mbox{a.s}\;\;
[P_{\theta_0}]\;\; \mbox{as}\;\; \ve \raro 0.
\ee
\vsp
\NI{\bf Proof :} This theorem follows by observing that the process
\be
R_\ve^N \equiv \sum_{i=1}^N \int_0^T \ve \sqrt{\lmd_i+1} Q_{i\ve}(t)dM_t^H,T \geq 0
\ee
is a local martingale with the quadratic variation process
\be
<R_\ve^N>_T=  \sum_{i=1}^N \int_0^T\ve^2 (\lmd_i+1)Q_{i \ve}^2(t)dw_t^H
\ee
and applying the Strong law of large numbers (cf. Liptser (1980); Prakasa
Rao (1999b), p. 61) under the condition (3.18) stated above.
\vsp

\NI {\bf Limiting distribution :}
\vsp

We now discuss the limiting distribution of the MLE $\hat \theta_{N \ve}$
as $\ve \raro 0.$
\vsp
\NI {\bf Theorem 3.2 :} Assume that the process $\{R_\ve^N, \ve \geq 0\}$ is a
local continuous martingale and that there exists a norming function
$I_\ve^N, \ve  \geq 0$ such that
\bea
\;\;\;\\ \nnb
(I_\ve^N)^2 <R_\ve^N>_T = (I_\ve^N)^2 \sum_{i=1}^N \int_0^T \ve^2 (\lmd_i+1)
Q_{i\ve}^2(t)dw_t^H \raro \eta^2 \;\; \mbox{in
probability}\;\; \mbox{as}\;\; \ve \raro 0 \\ \nnb
\eea
where $\eta$ is a random variable such that $P(\eta >0)=1.$ Then
\be
(I_\ve^N R_\ve^N, (I_\ve^N)^2 <R_\ve^N>_T) \raro (\eta Z, \eta^2) \mbox{ in
law}\;\; \mbox{as}\;\; \epsilon  \raro 0
\ee
where the random variable $Z$ has the standard Gaussian  distribution and the
random variables $Z$ and $\eta$ are independent.
\vsp
\NI {\bf Proof :} This theorem follows as a consequence of the central limit
theorem for local martingales (cf. Theorem 1.49 ; Remark 1.47 , Prakasa
Rao(1999b), p. 65).
\vsp
Observe that
\be
(I_\ve^N)^{-1}(\hat \theta_{N \ve} - \theta_0)  =  \frac{I_\ve^N R_\ve^N}
{(I_\ve^N)^2<R_\ve^N>}.
\ee
Applying the above theorem, we obtain the following result.
\vsp
\NI {\bf Theorem 3.3 :} Suppose the conditions stated in  Theorem 3.2
hold. Then
\be
(I_\ve^N)^{-1}(\hat \theta_{N \ve} - \theta_0) \raro \frac{ Z}{\eta} \mbox{ in
law}\;\; \mbox{as}\;\; \ve \raro 0
\ee
where the random variable $Z$ has the standard Gaussian  distribution and the
random variables $Z$ and $\eta$ are independent.
\vsp
\NI {\bf Remarks :} (i) If the random variable $\eta $ is a constant with probability one, then the
limiting distribution of the maximum likelihood estimator is
Gaussian  with mean 0 and variance $\eta^{-2}.$ Otherwise it is a mixture of
the Gaussian  distributions with mean zero and variance $\eta^{-2}$ with the
mixing distribution as that of $\eta.$
\vsp
(ii) Suppose that
\be
\lim_{N \raro \ity}\lim_{\ve \raro 0}\ve^2 I_\epsilon^N = I(\theta)
\ee
exists and is positive. Since the sequence of Radon-Nikodym derivatives
$$\{\frac{dP_{\tha}^{N,T,\ve}}{dP_{\tha_0}^{N,T,\ve}},n \geq 1 \}$$
form a non-negative martingale with respect to the filtration generated by the
sequence of random variables
$\{u_\ve^N, N \geq 1\}$, it converges almost surely to a
random variable $\nu_{\ve, \theta, \theta_0}$ as $N \raro \ity $  for every
$\ve >0.$ It is easy to see that the limiting random variable  is given by
\bea
\;\;\;\\\nonumber
\lefteqn{\nu_{\ve, \theta, \theta_0}(u_\ve)}\\\nonumber
&= & \exp \{ \sum_{i=1}^\ity \frac{\lmd_i+ 1}{\ve^2} [(\tha - \tha_0) \int_0^T Q_{i \ve}(t) dZ_{i\ve}(t)\\\nonumber  & &  \;\;\;\;-\frac{1}{2}\{(\tha - \lmd_i)^2 - (\tha_0 -\lmd_i)^2 \} \int_0^T
Q_{i\ve}^2(t) dw_t^H]\}.\\\nonumber
\eea
Furthermore  the sequence of random variables $u_{\ve}^N(t)$ converge in
probability to the random variable $u_{\ve}(t)$ as $N \raro \ity $ for every $\ve >0.$ Hence, by Lemma 4
in Skorokhod (1965, p. 100), it follows that the measures $P_{\theta}^\ve$ generated by the processes $u_\ve $for different values of $\theta,$ are absolutely continuous with respect to each other and the Radon-Nikodym
derivative of the probability measure $P_{\theta}^\ve $ with respect to the probability measure $P_{\theta_0}^\ve $
is given by
\bea
\frac{dP_{\theta}^\ve}{dP_{\theta_0}^\ve}(u_\ve) &=& \nu_{\ve, \theta,
\theta_0}(u_\ve) \\ \nnb
& = & \exp\{\sum_{i=1}^\ity \frac{\lmd_i + 1}{\ve^2} [(\tha - \tha_0)
\int_0^T Q_{i\ve}(t) dZ_{i\ve}(t)\\ \nnb
& & \;\; - \frac{1}{2} \{(\tha - \lmd_i)^2 - (\tha_0 - \lmd_i)^2 \}
\int_0^T Q_{i\ve}^2 (t) dw_t^H]\} .
\eea
It can be checked that the MLE $\hat{\tha}_{\ve}$ of $\tha $ based on
$u_{\ve}$ satisfies the likelihood equation
\be
\al_{\ve} = \ve^{-1} (\hat{\tha}_{\ve} - \tha_0) \beta_{\ve}
\ee
when $\tha_0$ is the true parameter where
\be
\al_{\ve} = \sum_{i=1}^{\ity} \sqrt{\lmd_i + 1} \int_0^T Q_{i\ve}(t) d
M_i^H(t)
\ee
and
\be
\bta_{\ve} = \sum_{i=1}^{\ity} (\lmd_i +1 ) \int_0^T Q^2_{i \ve}(t) dw_t^H .
\ee
\vsp
One can obtain sufficient conditions for studying the asymptotic behaviour
of the estimator $\hat{\tha}_{\ve}$ as in the finite projection case
discussed above. We omit the details.
\vsp

\NI {\bf Stochastic PDE with linear drift (singular case) : } \vsp
Let $(\Omg, \clf, P)$ be a probability space and consider the
process $u_{\ve} (t, x), 0 \leq x \leq 1, 0 \leq t \leq T$
governed by the stochastic partial differential equation \be d
u_{\ve}(t, x) = \tha \;  \triangle u_{\ve} (t,x) dt + \ve (I -
\triangle)^{-1/2} d \tilde W(t,x) \ee where $\tha > 0$ satisfying the
initial and the boundary conditions \bea u_{\ve}(0,x) &=& f(x),
\;0 < x <1, \; f \in L_2 [0,1] , \\ \nnb u_{\ve}(t,0) &=&
u_{\ve}(t,1) = 0, \;0 \leq t \leq T . \eea Here $I$ is the
identity operator, $\triangle = \frac{\pal^2}{\pal x^2}$ as
defined above and the process $\tilde W(t,x)$ is the cylindrical infinite
dimensional mfBm with $H \in [\frac{1}{2},1).$  Following the discussion in the
previous section, we assume the existence of a square integrable
solution $u_{\ve}(t,x)$ for the equation (3.31) subject to the
boundary conditions (3.32). Then the Fourier coefficients
$u_{i\ve}(t)$ of satisfy the stochastic differential equations \be
d u_{i \ve}(t) = - \tha \lmd_i u_{i \ve}(t) dt +
\frac{\ve}{\sqrt{\lmd_i + 1}} d \tilde W_i^H (t), \; 0 \leq t \leq T, \ee
with \be u_{i \ve}(0) = v _i, v _i = \int_0^1 f(x) e_i  (x) d x .
\ee \vsp Let $u_{\ve}^{(N)}(t,x)$ be the projection of  $u_{\ve}
(t,x)$ onto the subspace spanned by $\{e_1, \cdots, e_N \}$ in
$L_2 [0,1] .$ In other words \be u_{\ve}^{(N)} (t, x) =
\sum_{i=1}^{N} u_{i \ve}(t) e_i (x) . \ee Let $P_{\tha}^{(\ve, N)}
$ be the probability measure generated by $u_{\ve}^{(N)}$ on the
subspace spanned by $\{e_1, \cdots, e_N \}$ in $ L_2[0,1] .$ It
can be shown that the measures $\{P_{\tha}^{(\ve, N)}, \tha \in
\Tha\}$ form an equivalent family and \bea \;\;\;\;\\\nonumber
\lefteqn{ \log \frac{d P_{\tha}^{(\ve, N)}}{d P_{\tha_{0}}^{(\ve,
N)}}(u_{\ve}^{(N)})} \\ \nnb &=& - \frac{1}{\ve^2} \sum_{i=1}^{N}
\lmd_i (\lmd_i +1) [ (\tha - \tha_0) \int_0^T Q_{i \ve}(t) dZ_{i
\ve}(t) - \frac{1}{2} (\tha - \tha_0)^2 \lmd_i \int_0^T Q^2_{i
\ve} (t) dw_t^H ] . \eea It can be checked that the MLE
$\hat{\tha}_{\ve, N}$ of $\tha$ based on $u_{\ve}^{(N)}$ satisfies
the likelihood equation \be \al_{\ve, N} = -
\ve^{-1}(\hat{\tha}_{\ve, N} - \tha_0) \bta_{\ve, N} \ee when
$\tha_0$ is the true parameter where \be \al_{\ve , N} =
\sum_{i=1}^{N} \lmd_i \sqrt{\lmd_i +1} \int_0^T Q_{i \ve}(t) d M_i^H
(t) \ee and \be \bta_{\ve , N} = \sum_{i=1}^{N} (\lmd_i +1)
\lmd_i^2 \int _0^T Q_{i, \ve}^2 (t) dw_t^H . \ee Asymptotic
properties of these estimators can be investigated as in the
previous example. We do not go into the details as the arguments
are similar. 
\vsp
\NI{\bf Remarks :} One can study the local asymptotic mixed
normality (LAMN) of the family of probability measures generated by the
log-likelihood ratio processes by the standard arguments as in Prakasa
Rao (1999b) and hence investigate the asymptotic efficiency of the
MLE using Hajek-Lecam type bounds. 
\vsp
\newsection{ Parametric estimation  for stochastic parabolic equations driven by infinite dimensional mfBm}
\vsp 
We now extend some work of Cialenco et al. (2009)
dealing with problems of estimation in models more general than
those discussed in the previous section.  We introduce some
notation. 
\vsp 
Let ${\bf H}$ be a separable Hilbert space with the
inner product $(.,.)_0$ and with the corresponding norm $||.||_0.$
Let $\Lambda $ be a densely defined linear operator on ${\bf H}$
with the property that there exists $c>0$ such that
$$||\Lambda u||_0\geq c||u||_0$$
for every $u$ in the domain of the operator $\Lambda.$ The operator powers $\Lambda ^{\gamma} ,\gamma \in R$
are well defined and generate the spaces ${\bf H}^\gamma$ with the properties (i) for $\gamma >0, {\bf H}^\gamma $ is
the domain of $\Lambda ^\gamma,$ (ii) ${\bf H}^0={\bf H}$, and (iii) for $\gamma <0, {\bf H}^\gamma $ is the
completion of ${\bf H}$ with respect to the norm $||.||_\gamma \equiv ||\Lambda^ \gamma .||_0$ (cf. Krein et al. (1982)).
The family of spaces $\{{\bf H}^\gamma, \gamma \in R\}$ has the following properties:

\begin{description}
\item {(i)} $\Lambda^{\gamma}({\bf H}^r)= {\bf H}^{r-\gamma}, \gamma, r \in R;$
\item{(ii)} For $\gamma_1 < \gamma_2,$ the space ${\bf H}^{\gamma_2}$ is densely and continuously embedded
into  ${\bf H}^{\gamma_1}$, that is,  ${\bf H}^{\gamma_2} \subset{\bf H}^{\gamma_1}$ and there exists a
constant $c_{12} >0$ such that $||u||_{\gamma_1} \leq c_{12}||u||_{\gamma_2};$
\item{(iii)} for every $\gamma \in R$ and $m >0,$ the space  ${\bf H}^{\gamma -m}$ is the dual of the
space ${\bf H}^{\gamma +m}$ with respect to the inner product in ${\bf H}^{\gamma},$
with duality $<.,.>_{\gamma,m}$ given by
$$<u_1,u_2>_{\gamma,m}= (\Lambda^{\gamma -m}u_1, \Lambda^{\gamma + m}u_2)_0, \;\;u_1 \in {\bf H}^{\gamma -m}, u_2 \in {\bf H}^{\gamma +m}.$$
\end{description}
\vsp
Let $(\Omega, \clf,P)$ be a probability space and let $\{\tilde W_j^H, j \geq 1\}$ be a family of independent
mixed fractional Brownian motions on this space with the {\it same} Hurst index $H$ in $(0,1).$
\vsp
Consider the SDE
\be
du(t)+(\cla_0+\theta \cla_1)u(t)dt= \sum_{j\geq 1}g_j(t)d\tilde W_j^H(t), 0 \leq t \leq T, u(0)=u_0 
\ee
where $\cla_0, \cla_1$   are linear operators, $g_j, j \geq 1$ are non-random and $\theta \in \Theta\subset R.$
The equation (4.1) is said to be  {\it diagonalizable} if the operators $\cla_0, \cla_1$ have the same  system of
eigenfunctions $\{h_j, j \geq 1\}$ such that $\{h_j, j \geq 1\}$ is an orthonormal basis in ${\bf H}$
and each $h_j$ belongs to $\cap_{\gamma \in R}{\bf H}^\gamma.$ It is called $(m, \gamma)$-parabolic  for some $m \geq 0, \gamma \in R,$ if
\vsp
(i) the operator $\cla_0+\theta \cla_1$ is uniformly bounded from ${\bf H}^{\gamma+m}$ to ${\bf H}^{\gamma-m}$ for every $\theta \in \Theta $, that is, there exists $C_1 >0$ such that
\be
||(\cla_0+\theta \cla_1)v||_{\gamma-m}\leq C_1||v||_{\gamma+m},\;\; \theta \in \Theta, v \in H^{\gamma+m};
\ee
and\\
(ii)there exists a $\delta >0$ and $C \in R$ such that, 
\be
-2<(\cla_0+\theta \cla_1)v,v>_{\gamma,m}+\delta||v||^2_{\gamma+m} \leq C ||v||^2_{\gamma}, v \in {\bf H}^{\gamma+m}, \theta \in \Theta.
\ee
\vsp
If the equation (4.1) is $(m, \gamma)$-parabolic, then the condition (ii) implies that
$$<(2\cla_0+2\theta\cla_1+CI)v,v)_{\gamma,m} \geq \delta||v||^2_{\gamma+m}$$
where $I$ is the identity operator. The Cauchy-Schwartz inequality and the continuous embedding of
${\bf H}^{\gamma+m}$ into ${\bf H}^\gamma $ will imply that
$$ ||(2\cla_0+2\theta\cla_1+CI)v||_\gamma \geq \delta_1||v||_{\gamma}$$
for some $\delta_1 >0$ uniformly in $\theta \in \Theta.$ 
\vsp 
Let us choose $\Lambda=  [2\cla_0+2\theta_0 \cla_1+CI]^{1/(2m)}$ for
some fixed $\theta_0 \in \Theta.$ If the operator $\cla_0+\theta
\cla_1$ is unbounded, we say that $\cla_0+\theta \cla_1$ has order
$2m$ and $\Lambda $ has order 1. If the equation (4.1) is
$(m,\gamma)$-parabolic and diagonalizable, we will assume that the
operator $\Lambda $ has the same eigenfunctions as the operators
$\cla_0$ and $\cla_1.$ This is justified by the comments made
above. 
\vsp 
Suppose the equation (4.1) is diagonalizable and
there exists eigenvalues $\{\rho_j, j \geq 1\},\{\nu_j, j\geq 1\}$
such that
$$\cla_0 h_j=\rho_j h_j \;\;\mbox{and}\;\; \cla_1 h_j=\nu_j h_j.$$
Without loss of generality, we can also assume that there exists $\{\lambda_j, j \geq 1\}$ such that
$$\Lambda h_j= \lambda_j h_j. $$
Following the arguments in Cialenco et al. (2009), it can be shown that the equation (4.1) is $(m, \gamma)$-parabolic
if and only if there exists  $\delta >0, C_1 >0$ and $C_2 \in R$ such that, for all $j\geq 1, \theta, \in \Theta,$
\be
|\rho_j+\theta \nu_j|\leq C_1 \lambda_j^{2m} 
\ee
and
\be
-2(\rho_j+\theta \nu_j)+\delta \lambda_j^{2m} \leq C_2. 
\ee

As the conditions in (4.4) and  (4.5) do not depend on $\gamma,$
we conclude that a diagonalizable equation (4.1) is $(m,
\gamma)$-parabolic for some $\gamma $ if and only if it is
$(m,\gamma)$-parabolic for every $\gamma.$ Here after we will say
that the equation (4.1) is $m$-parabolic. We will assume that the
equation (4.1) is diagonalizable and fix the basis $\{h_j, j \geq
1\}$ in ${\bf H}$  consisting of the eigenfunctions of $\cla_0$,
$\cla_1 $ and $\lambda $. Recall that set of eigenfunctions is the
same for all the three operators. Since $h_j$ belongs to every
${\bf H}^{\gamma},$ and since $\cap_\gamma{\bf H}^{\gamma}$ is
dense in $\cup_\gamma{\bf H}^{\gamma},$ every element $f$ of
$\cup_\gamma{\bf H}^{\gamma},$ has a unique expansion $\sum_{j\geq
1}f_jh_j$ where $f_j=<f,h_j>_{0,m}$ for suitable $m.$ 

The functional structure described above follows the work in Cialenko et al. (2009).

\vsp

\NI{\bf Definition :} The {\it infinite dimensional mixed fractional Brownian motion} $\tilde W^H$ is
an element of $\cup_{\gamma \in R}{\bf H}^{\gamma}$ with the
expansion 
\be \tilde W^H(t)= \sum_{j\geq 1}\tilde h_j \tilde W_j^H(t). 
\ee 
 \NI{\bf Definition :} The solution of the diagonalizable equation 
\be
du(t)+(\cla_0+\theta \cla_1)u(t)dt=d\tilde W^H(t), 0 \leq t \leq T,
u(0)=u_0, 
\ee 
with $u_0\in {\bf H},$ is defined to be a random process $u(t), 0< t \leq T,$ with
values in $\cup_\gamma{\bf H}^{\gamma}$ and has an expansion 
\be
u(t)= \sum_{j\geq 1}h_j u_j(t) 
\ee 
where 
\be u_j(t)=(u_0,h_j)_0 e^{-(\theta \nu_j+\rho_j)t}+ \int_0^t e^{-(\theta \nu_j+\rho_j)(t-s)}\;\;d \tilde W_j^H(s). 
\ee
Let 
\be 
\mu_j(\theta)=\theta \nu_j+\rho_j, j \geq 1. 
\ee 
In view of (4.5), we get that there exists a positive integer $ J$ such
that 
\be 
\mu_j(\theta)>0 \;\;\mbox{for}\;\; j \geq J 
\ee 
if the equation (4.1) is $m$-parabolic and diagonalizable. 
\vsp 
\NI{\bf Theorem 4.1 :} {\it Suppose that $H \geq \frac{1}{2}$ and the equation
(4.1) is $m$-parabolic and diagonalizable. Further suppose that
there exists a positive real number $\gamma$ such that} 
\be 
\sum_{j \geq 1}(1+|\mu_j(\theta)|)^{-\gamma} < \ity.  
\ee 
{\it Then, for every} $t >0, \tilde W^H(t) \in L_2(\Omega, {\bf H}^{-m \gamma}) \;\mbox{and}\; u(t) \in L_2(\Omega, {\bf H}^{-m \gamma + 
m \min(2H,1)}).$
\vsp 
\NI{\bf Proof :} The condition (3) implies that $\lim_{j \raro \ity}|\mu_j|=\ity,$ and hence the operators $\cla_0+\theta \cla_1$ and $\Lambda$ are unbounded. The parabolicity assumption and the equations (4.4) and (4.5) imply that, for sufficently large $j,$
$$1+|\mu_j(\theta)|^m\leq C_2\lambda_j^{2m}$$
uniformly in $\theta\in \Theta.$ Furthermore
$$E||\tilde W^H(t)||^2_{-m\gamma}= (t^{2H}+t)\sum_{j\geq 1}\lambda_j^{-2m\gamma} \leq C_2(t^{2H}+t)\sum_{j\geq 1}(1+\mu_j(\theta)|^{-\gamma}<\ity.$$ From the definition of the mixed fractional Brownian motion and the fact that the component fractional Brownian motion and the Brownian motion are independent and centered, it follows that
\bea
E[u_j^2(t)]&= & H(2H-1)e^{-2\mu_j(\theta)t}\int_0^t \int_0^t e^{\mu_j(\theta)(s_1+s_2)}|s_1-s_2|^{2H-2}ds_1 ds_2 + e^{-2\mu_j(\theta)t}\int_0^te^{2\mu_j(s)}ds\\\nonumber
&= & J_1(t)+ J_2(t) \;\;\mbox{(say).}\\\nonumber
\eea
by the properties of fractional Brownian motion (cf. Prakasa Rao (2010)) and Brownian motion. It can be checked that
$$\lim_{j \raro \ity}|\mu_j(\theta)|^{2H}J_1(t)= H(2H-1)\int_0^\ity x^{2H-2}e^{-x}dx=H(2H-1)\Gamma(2H-1)$$
and
$$\lim_{j \raro \ity}|\mu_j(\theta)|J_2(t)= \frac{1}{2},$$
Hence
$$\lim_{j\raro \ity} |\mu_j(\theta)|^{\min(2H,1)}E(u_j^2(t)) <\ity.$$
and
$$\sum_{j=1}^\ity (1+|\mu_j(\theta)|)^{-\gamma+\min(2H,1)}E(u_j^2(t))<\ity. $$
\vsp
\NI{\bf Maximum likelihood estimation :} 
\vsp 
Consider the diagonalizable equation 
\be
du(t)+(\cla_0+\theta \cla_1)u(t)dt= d\tilde W^H(t), u(0)=0, 0 \leq t \leq T 
\ee 
with 
\be u(t) =\sum_{j\geq 1}h_j u_j(t)
\ee 
as given by (4.9). Suppose the processes $\{u_i(t), 0 \leq t \leq T\}, i=1,\dots,N $ can be
observed continuously over the interval $[0,T].$ The problem is to estimate the
parameter $\theta$ using these observed paths over the interbal $[0,T].$. \vsp Note that
$\mu_j(\theta)= \rho_j+\nu_j\theta,$ where $\rho_j$ and $\nu_j$ are
the eigenvalues of $\cla_0$ and $\cla_1$ respectively. Furthermore
each process $u_j$ is a mixed fractional Ornstein-Uhlenbeck process
satisfying the stochastic differential equation 
\be
du_j(t)=-\mu_j(\theta)u_j(t)dt+d \tilde W_j^H(t),u_j(0)=0, 0 \leq t \leq T. 
\ee 
Since the processes $\{\tilde W_j^H, j \geq 1\}$ are independent,
it follows that the processes $\{u_j, 1 \leq j \leq N\}$ are
independent. Following the notation introduced above, let
\be 
M_j^H(t)= \int_0^tg_H(s,t)d\tilde W_j^H(s), Q_j(t)=\frac{d}{dw_j^H(t)}\int_0^tg_H(s,t)u_j(s)ds 
\ee 
and 
\be 
Z_j(t)=\int_0^tg_H(s,t)du_j(s) 
\ee for $j=1,\dots,N.$ Applying the Girsanov-type formula, it can be shown that
the measure generated by the process $(u_1, \dots,u_N)$ is
absolutely continuous with respect to the measure generated by the
process $(\tilde W_1^H,\dots,\tilde W_N^H)$ and their  Radon-Nikodym derivative
is given by 
\be 
\exp(-\sum_{j=1}^N \mu_j(\theta)\int_0^T Q_j(s)dZ_j(s)-\sum_{j=1}^N \frac{[\mu_j(\theta)]^2}{2}\int_0^TQ_j^2(s)dw_H(s)). 
\ee
Maximizing this function with respect to the parameter $\theta,$ we get the maximum likelihood estimator
$$\hat \theta_N= -\frac{ \sum_{j=1}^N\int_0^T \nu_j Q_j(s)(dZ_j(s)+\rho_jQ_j(s)dw_H(s))}{\sum_{j=1}^N\int_0^T\nu_j^2 Q_j^2(s)dw_H(s)}.$$
\vsp
\NI{\bf Theorem 4.2 :} {\it Suppose the Hurst index $ H>\frac{1}{2}.$ Suppose that} 
\be
\sum_{j=1}^\ity \frac{\nu_j^2}{\mu_j(\theta)}= \ity. 
\ee
Then
$$\lim_{N \raro \ity}\hat \theta_N=\theta \;\;\mbox{a.s.}$$

{\it where $\theta $ is the true parameter. Furthermore, if the
condition (i) holds, then} 
\be 
\lim_{N \raro \ity}(\sum_{j=J}^N\frac{\nu_j^2}{\mu_j(\theta)})^{1/2}(\hat \theta_N-\theta)
\stackrel \cll \raro \;\;N(0,1) \mbox{as}\;\; N \raro \ity 
\ee 
{\it where} $J=\min\{j:\mu_i(\theta)=0 \;\;\mbox{for all}\;\;i \geq j\}.$ 
\vsp 
\NI{\bf Proof :} Observe that 
\be
\hat \theta_N-\theta= -\frac{\sum_{j=1}^N\int_0^T\nu_jQ_j(s)dM_j^H(s)}{\sum_{j=1}^N\int_0^T\nu_j^2Q_j^2(s)dw^H(s)}. 
\ee
The numerator of the expression on the right side (4.22) is a local martingale and the denominator is its quadratic variation. It is known that  
\be 
\lim_{j \raro \ity}\mu_j(\theta)\;\;E[\int_0^TQ_j^2(s)dw_H(s)]= \frac{T}{2} >0 
\ee
from the equation (2.1)  in Chigansky and Kleptsyna (2019) (cf. Maruskeevych (2016)) which implies that   

$$\sum_{j=1}^N\int_0^T\nu_j^2Q_j^2(s)dw^H(s)\raro \ity \;\;\;\mbox{a.s. as} \;\;N \raro \ity$$ 
under the condition (4.20). Hence, an application of the Strong law of large numbers for local martingales (cf. Prakasa Rao (1999), Liptser and Shiryayev (1989)) will imply that the left side of the equation (4.22) converges to zero almost surely and hence
$$\hat \theta_N-\theta \raro 0 \;\;\mbox{almost surely}$$
as $N \raro \ity.$
\vsp
\NI{\bf Theorem 4.3:} {\it Suppose there exists a norming function $I_N$ such that 
\be
I_N^2 \sum_{j=1}^N\int_0^T\nu_j^2Q_j^2(s)dw^H(s) \raro \eta^2 \;\;\mbox{in probability as} \;\;N\raro \ity
\ee
where $\eta$ is a random variable such that $P(\eta>0)=1.$ Let 
\be
R_N= \sum_{j=1}^N\int_0^T\nu_jQ_j(s)dM_j^H(s).
\ee
Then
\be
(I_N R_n, I_N^2 <R_N>)\raro (\eta Z,\eta^2) \;\;\mbox{in law as }\;\;N \raro \ity
\ee
where the random variable $Z$ has the standard normal distribution and the random vatiables $Z$ and $\eta$ are independent.Here $<R_N>, N\geq 1$ is the quadratic variation of the process $\{R_N, N\geq 1\}.$}
\vsp
\NI{\bf Proof:} This theorem follows as a consequence of the central limit theorem for local martingales (cf. Theorem 1.49 and remark 1.47, Prakasa Rao (1999b)).
\vsp
Observe that
\be
I_N^{-1}(\hat \theta_N-\theta)= \frac{I_N R_N}{I_N^2 <R>_N}, N \geq 1.
\ee
Applying Theorem 4.3, we obtain the following result .
\vsp
\NI{\bf Theorem 4.4:} {\it Suppose the conditions stated in Theorem 4.2 hold. Then
\be
I_N^{-1}(\hat \theta_N-\theta) \raro \frac{Z}{\eta} \;\;\mbox{in law as} \;\;N \raro \ity
\ee
where the random variable $Z$ has the standard normal distribution and the random variable $Z$ and $\eta$ are independent.}
\vsp
\NI{\bf Remarks :} If the random variable $\eta$ is a constant almost surely, then the limiting distribution of the maximum likelihood estimator $\hat \theta_N$ is Gaussian with mean zero and variance $\eta^{-2}.$ Otherwise it is a mixture of the normal distribution with mean zero and variance $\eta^{-2}$ with the mixing distribution as that the random variable $\eta.$ Applying Berry-Esseen bound for sums of independent random variables, it is possible to obtain the limiting distribution and rate of convergence of the distribution of the MLE if a bound on the variance of the random variable
$$\int_0^TQ_j^2(s)\;dw^H(s)$$
can be obtained as $j \raro \ity$ after suitable norming as in Lemma A.2 in Cialenco et al. (2009) in the fractional Brownian case. It has not been possible to obtain a similar result in the case mixed fractional Brownian motion as the function $g_H(s,t)$ defining the martingale $M_j^H$ defined by (3.8) is not explicitly known. The problem of obtaining a closed form for the function $g_H(s,t)$ remains open.
\vsp
\NI{\bf Acknowledgment} This work was supported under the scheme ``INSA Senior Scientist" by the Indian National Science Academy (INSA) at the CR Rao Advanced Institute of Mathematics, Statistics and Computer Science, Hyderabad, India.
\vsp
\NI{\bf References :}
\begin{description}
\item Cai, C., Chigansky., and Kleptsyna, M. (2016)  Mixed Gaussian processes ; A filtering approach, {\it Ann. Probab.}, {\bf 44}, 3032-3075.

\item Cheridito, C. (2001) Mixed fractional Brownian motion, {\it Bernoulli}, {\bf 7}, 913-934.

\item Chigansky, P., and Kleptsyna, M. (2019) Statistical analysis of the mixed fractional Ornstein-Uhlenbeck process. {\it Theory Probab. Appl.}, {\bf 63}, 408-425.

\item Cialenco, I. (2018) Statistical inference for SPDEs: an overview, {\it Statist. Infer. Stoch. Proc.}, https://doi.org/10.1007/s11203-018-9177-9.

\item Cialenco, I., Lototsky, S.G., and Pospisil, J. (2009) Asymptotic properties of the maximum likelihood estimator for
stochastic parabolic equations with additive fractional Brownian motion, {\it Stoch. and Dynam.}, {\bf 9}, 169-185.

\item da Silva, Jose Luis., Erroui, M and essaky, El Hassan (2018) Mixed stochastic differential equations: Existence and uniqueness result, {\it J. Theor. Probab.}, {\bf 31}, 1119-1141.

\item Da Prato, G. and Zabczyk, J. (1992) {\it Stochastic Equations in
Infinite Dimensions}, Cambridge University Press.

\item Guerra, Joao., and Nualart, D. (2008) Stochastic differential equations driven by fractional Brownian motion and standard Brownian motion, {\it Stoch. Anal Appl.}, {\bf 26}, 1053-1075.

\item Huebner, M., Khasminski, R. and Rozovskii. B.L. (1993) Two
examples of parameter estimation for stochastic partial differential
equations, In {\it Stochastic Processes : A Festschrift in Honour of
Gopinath Kallianpur}, Springer, New York, pp. 149-160.

\item Huebner, M., and Rozovskii, B.L. (1995) On asymptotic properties of
maximum likelihood estimators for parabolic stochastic SPDE's, {\it  Prob.
Theory and Relat. Fields},{\bf  103}, 143-163.

\item Ito, K. (1984) Foundations of Stochastic Differential Equations
in Infinite Dimensional Spaces, Vol. {\bf 47}, CBMS Notes, SIAM, Baton
Rouge.

\item Kallianpur, G., and Xiong, J. (1995) {\it Stochastic Differential
Equations in Infinite Dimensions} ,  IMS Lecture Notes, Vol.{\bf 26}, Hayward,
California.

\item Kleptsyna, M. L. and Le Breton, A. (2002) Statistical analysis of
the fractional Ornstein-Uhlenbeck type process, {\it Statist.
Infer. for Stoch. Proc.}, {\bf 5}, 229--248.

\item Krein, S.G., Petunin, Yu.I., and Semenov, E.M. (1982)
{\it Interpolation of linear operators}, Volume 54 of {\it Translations of Mathematical Monographs},
American Mathematical Society, Providence, Rhode Island.

\item Le Breton, A. (1998) Filtering and parameter estimation in a simple
linear model driven by a fractional Brownian motion, {\it Statist. Probab.
Lett.},{\bf 38}, 263-274.

\item Liptser, R. (1980) A strong law of large numbers, {\it Stochastics},
{\bf 3}, 217-228.

\item Lototsky, S.V. and Rozovsky, B.L. (2017) {\it Stochastic Partial Differential Equations}, Springer, Switzerland.

\item Marushkevych, Dmytro (2016) Large deviations for drift parameter estimator of mixed fractional Ornstein-uhlenbeck process, {\it Modern Stochastics: Theory and applications}, {\bf 3}, 107-117.

\item Mishura, Y., Ralchenko, K., and Shevchenko, G. (2019) Existence and uniqueness of mild solutions to the stochastic heat equation with white and fractional noises, {\it Theory of Prob. and Math. stat.}, {\bf 98} 149-170.

\item Mishura, Y., and Shevchenko, G. (2011) Existence and uniqueness of the solution of stochastic differential equation involving wiener process and fractional Brownian motion with hurst index $H>\frac{1}{2}.$  {\it Commun. Stat. Theory Methods}, {\bf 40}, 3492-3508.

\item  Prakasa Rao, B.L.S. (1987) {\it Asymptotic Theory of Statistical
Inference}, Wiley, New York.

\item  Prakasa Rao, B.L.S. (1999a) {\it Statistical Inference for
Diffusion Type Processes}, Arnold, London and Oxford University
Press, New York.

\item  Prakasa Rao, B.L.S. (1999b) {\it Semimartingales and Their
Statistical Inference}, CRC Press, Boca Raton and Chapman and
Hall, London.

\item Prakasa Rao, B.L.S. (2000) Bayes estimation for stochastic partial
differential equations, {\it J. Statist. Plan. Inf.},{\bf 91}, 511-524.

\item Prakasa Rao, B.L.S. (2001) Statistical inference for stochastic
partial differential equations, In {\it Selected Proceedings of the
Symposium on Inference for Stochastic Processes}, Ed. I.V.Basawa,
C.C.Heyde and R.L.Taylor, IMS Monograph Series,  Vol. {\bf 37}, pp. 47-70.

\item Prakasa Rao, B.L.S. (2002a) On some problems of estimation for some
stochastic partial differential equations, In {\it Uncertainty and
Optimality}, Ed. J.C.Mishra, World Scientific, Singapore, pp. 71-153.

\item Prakasa Rao, B.L.S. (2002b) Minimum distance estimation for some stochastic partial differential equations, {\it J. Korean Stat. Soc.} {\bf 31} 213-228.

\item Prakasa Rao, B.L.S. (2003) Parameter estimation for linear stochastic differential equations driven by fractional Brownian motion, {\it Random Operators and Stochastic Equations}, {\bf 11}, 229-242.

\item Prakasa Rao, B.L.S. (2004) Parameter estimation for some stochastic partial differential
equations driven by infinite dimensional fractional Brownian motion, {\it Theory Stochastic. Process}, {\bf 10 (26)} 116-125.

\item Prakasa Rao, B.L.S. (2009) Estimation for stochastic differential equations driven by mixed fractional Brownian motion. {\it Calcutta Stat. Assoc. Bull.} 61: 143-153.

\item Prakasa Rao, B.L.S. (2010) {\it Statistical Inference for Fractional Diffusion Processes}, Wiley, London.

\item Prakasa Rao, B.L.S.  (2013) Parameter estimation for a two-dimensional  stochastic Navier-Stokes equation driven by infinite dimensional fractional Brownian motion, {\it Random Operators and Stochastic Equations}, {\bf 21}  37-52.

\item Prakasa Rao, B.L.S. (2015a) Option pricing for processes driven by mixed fractional Brownian motion with superimposed jumps. {\it Probability in the Engineering and Information Sciences}. 29:  589-596.

\item Prakasa Rao, B.L.S. (2015b) Pricing geometric Asian power options under mixed fractional Brownian motion environment, {\it Physica A}, {\bf 446}, 92-99.

\item Prakasa Rao, B.L.S. (2017a) Instrumental variable estimation for a linear stochastic differential equation driven by a mixed fractional Brownian motion. {\it Stochastic Anal. Appl.} 35: 943-953.

\item Prakasa Rao, B.L.S. (2017b) Optimal estimation of a signal perturbed by a mixed fractional Brownian motion, {\it Theory of Stochastic Processes}, {\bf 22 (38)}, 62-68.

\item Prakasa Rao, B.L.S. (2018a) Parametric estimation for linear stochastic differential equations driven by mixed fractional Brownian motion, {\it Stochastic Analysis and Applications}, {\bf 36}, 767-781.

\item Prakasa Rao, B.L.S. (2018b) Pricing geometric Asian options under mixed fractional Brownian motion environment with superimposed jumps, {\it Calcutta Statistical Association Bulletin}, {\bf 70}, 1-6.

\item Prakasa Rao, B.L.S. (2019) Nonparametric estimation of trend for stochastic differential equations driven by mixed fractional Brownian motion, {\it Stochastic Analysis and Applications}, {\bf 37}, 271-280.

\item Prakasa Rao, B.L.S. (2020) Nonparametric estimation for stochastic differential equations driven by mixed fractional Brownian motion  with random effects,  In the Special Issue in honour of CR Rao Birth Centenary, {\it Sankhya, Series A} (to appear).

\item Prakasa Rao, B.L.S. (2021a) Maximum likelihood estimation in the mixed fractional Vasicek model, {\it Journal of Indian Society for Probability and Statistics}.

\item Prakasa Rao, B.L.S. (2021b) Nonparametric estimation of linear multiplier for processes driven by mixed fractional Brownian motion,  In the Special Issue in memory of Aloke Dey, {\it Statistics and Applications}, {\bf 19} 1-12 (to appear).

\item Rozovskii, B.L. (1990) {\it Stochastic Evolution Systems}, Kluwer, Dordrecht.

\item Samko, S.G., Kilbas, A.A., and Marichev, O.I. (1993) {\it Fractional
Integrals and derivatives}, Gordon and Breach Science.

\item Skorokhod, A.V. (1965) {\it Studies in the Theory of Random Processes}, Addison-Wesley, Reading, MA. 

\end{description}
\end{document}